\documentclass[a4paper,leqno]{amsart}
       
\usepackage{amsmath,amsfonts,amssymb}
\usepackage{amsthm}
\usepackage{hyperref}

\def\C{{\mathbb C}}
\def\R{{\mathbb R}}
\def\N{{\mathbb N}}

\def\Sph{{\mathbb S}}
 
\def\virgp{\raise 2pt\hbox{,}}

\def\({\left(}
\def\){\right)}
\def\<{\left\langle}
\def\>{\right\rangle}
\def\le{\leqslant}
\def\ge{\geqslant}

\def\d{{\partial}}

\def\om{\omega}
\def\si{{\sigma}}

\def\F{\mathcal F}

\def\O{\mathcal O}

\def\ii{\int\!\!\!\int}
\DeclareMathOperator{\RE}{Re}
\DeclareMathOperator{\IM}{Im}

\theoremstyle{plain}
\newtheorem{theorem}{Theorem}[section]
\newtheorem{lemma}[theorem]{Lemma}

\newtheorem{proposition}[theorem]{Proposition}
\newtheorem*{hyp}{Assumption}

\theoremstyle{definition}
\newtheorem*{definition}{Definition}

\theoremstyle{remark}
\newtheorem{remark}[theorem]{Remark}
\newtheorem*{remark*}{Remark}
\newtheorem*{example}{Example}

\numberwithin{equation}{section}

\begin{document}

\title[Cauchy problem in $H^s$ for NLS
  with potential]{On the Cauchy problem in Sobolev spaces for
  nonlinear Schr\"odinger equations with potential}  
\author[R. Carles]{R{\'e}mi Carles}
\address{Wolfgang Pauli Institute, c/o Inst.~f.~Math.\\ Universit\"at Wien\\
        Nordbergstr.~15\\ A-1090 Wien\\ Austria\footnote{Future
  affiliation: Univ. Montpellier,
  UMR CNRS 5149,  34095
  Montpellier cedex 5, France.}}
\email{Remi.Carles@math.cnrs.fr}
\thanks{Support by the French ANR
  project SCASEN is acknowledged.} 
\begin{abstract}
We consider the Cauchy problem for nonlinear Schr\"odinger equations
in the presence of a smooth, possibly unbounded, potential. No
assumption is made on the sign of the potential. If the potential
grows at most linearly at infinity, we construct solutions in Sobolev
spaces (without weight), locally in time. Under some natural
assumptions, we prove that the $H^1$-solutions are global in time. On
the other hand, if the potential has a super-linear growth, then the
Sobolev regularity of positive order is lost instantly, not matter how
large it is, 
unless the initial datum decays sufficiently fast at infinity.   
\end{abstract}
\maketitle

\section{Introduction}
\label{sec:intro}
We consider the Cauchy problem for the (nonlinear) Schr\"odinger
equation
\begin{equation}
  \label{eq:nlsp}
  \begin{aligned}
    i\d_t u +\frac{1}{2}\Delta u &= V(x) u + f\( |u|^2\) u,\quad
    (t,x)\in \R_+\times \R^d,\\
    u_{\mid t=0}&= a_0 \in H^s\(\R^d\),\quad s\ge 0,
  \end{aligned}
\end{equation}
where the potential $V$ is smooth and sub-quadratic (see below), the
nonlinearity $f$ is sufficiently smooth, and the initial data $a_0$
may or may not belong to weighted $L^2$ spaces $\F (H^k)$ (sometimes
denoted $L^2_k$), where
$\F$ stands for the Fourier transform. Note that we consider only
propagation in the future; this choice is made only to simplify some
statements. We show that if the potential
$V$ is sub-linear, then \eqref{eq:nlsp} is locally well-posed in $H^1(\R^d)$,
upon suitable assumptions on $f$. On the other hand, if $V$ is
super-linear (e.g. harmonic potential), then \eqref{eq:nlsp} is
ill-posed in all Sobolev spaces of positive order; this is not a
nonlinear result, since 
it holds even when $f\equiv 0$. This is heuristically reasonable, at
least in the case of the harmonic oscillator: the potential rotates
the phase space, so the natural space for the initial data is of the
form $H^s\cap \F(H^s)$. If $a_0\in H^s\setminus \F(H^k)$, for $s\ge
k>0$, then $u(t,\cdot )\not \in H^k(\R^d)$ for arbitrarily 
small $t>0$. 
For the linear equation, this can be seen via the Fourier
integral representation (Mehler's formula in the case of the harmonic
potential). The proof we present treats both linear and nonlinear
cases.
\smallbreak

Before going further into details, we clarify our assumptions. We
define the Fourier transform as
\begin{equation*}
  \F \varphi(\xi) =\widehat \varphi(\xi)= \int_{\R^d} e^{-ix\cdot
  \xi}\varphi(x)dx. 
\end{equation*}
Denote $\<x\> =(1+|x|^2)^{1/2}$. For $s\ge 0$, we define
\begin{equation*}
  H^s\(\R^d\)=\{ \varphi\in L^2(\R^d)\ ;\ \xi\mapsto \<\xi\>^s \widehat
  \varphi(\xi) \in L^2(\R^d)\} \quad ;\quad H^\infty \(\R^d\)=\cap_{s\ge
  0}H^s\(\R^d\). 
\end{equation*}
In particular,  $\F(H^s)$ is just the weighted $L^2$ space:
\begin{equation*}
 \F(H^s)=  \{ \varphi\in L^2(\R^d)\ ;\ x\mapsto \<x\>^s 
  \varphi(x) \in L^2(\R^d)\}.
\end{equation*}
\begin{hyp}
  We assume that the potential is smooth, real-valued and
  sub-quadratic: 
  $V\in C^\infty(\R^d;\R)$, and $\partial^\alpha V\in
  L^\infty(\R^d)$ for all $|\alpha|\ge 2$.
\end{hyp}
\begin{definition}
  We say that $V$ is \emph{sub-linear} if $\partial^\alpha V\in
  L^\infty(\R^d)$ as soon as $|\alpha|\ge
  1$. We say that $V$ is \emph{super-linear} if
  $\nabla_x V$ is unbounded. 
\end{definition}
\begin{remark*}
For super-quadratic potentials, the
theory must be modified. First, if $V$ is super-quadratic and
negative, then $H$ is not essentially self-adjoint on
$C_0^\infty(\R^d)$ (\cite{Dunford,ReedSimon2}). If  $V$ is super-quadratic and
positive, then even the local existence results are different. We
refer to \cite{YajZha01,YajZha04} for very interesting results in this
direction.  
\end{remark*}
 
The construction of the parametrix for the propagator
of $H=-\frac{1}{2}\Delta +V$ provided by D.~Fujiwara
\cite{Fujiwara79,Fujiwara} shows that $U(t)=e^{-itH}$, which is
$L^2$-unitary, satisfies the local dispersion estimate: there exists
$\delta>0$ such that
\begin{equation}\label{eq:dispfujiwara}
  \| U(t)\|_{L^1\to L^\infty}\lesssim \frac{1}{|t|^{d/2}},\text{ for
  }|t|\le \delta.
\end{equation}
One can infer local and global existence results for
\eqref{eq:nlsp} if $a_0\in D( \sqrt H)$ when $V\ge 0$, under
suitable assumptions on the nonlinearity $f$, as proved initially by
Oh \cite{Oh}. The assumption $V\ge 0$ is actually not necessary, and
one can prove the local existence results of Oh in weighted Sobolev
spaces  of the form $H^s\cap \F(H^s)$ 
thanks to Strichartz estimates (see
e.g. \cite{CazCourant}, and \cite{CaCCM} where global existence
results are recalled for potentials $V$ which are not necessarily
non-negative). In all this paper, $u$ is assumed to be a mild solution
to \eqref{eq:nlsp}, that is, to solve
\begin{equation*}
  u(t) = U(t)a_0 - i\int_0^t U(t-s)\( f\(|u(s)|^2\)u(s)\)ds. 
\end{equation*}
In Proposition~\ref{prop:local} though, we construct a classical
solution for \eqref{eq:nlsp}. 

\smallbreak
When $V\ge 0$ and $f(|u|^2)=\mu 
|u|^{2\si}$, one can prove global 
existence in $D(\sqrt H)$ for the solution $u$ of \eqref{eq:nlsp}
under suitable assumptions on $\mu$ and $\si$,
thanks to the following conservations:
\begin{align*}
  \text{Mass: }&\frac{d}{dt}\(\| u(t)\|_{L^2}^2\)= 0.\\
\text{Energy: }&\frac{d}{dt}\(\frac{1}{2}\|\nabla
u(t)\|_{L^2}^2+\frac{\mu}{\si+1} \|u(t)\|_{L^{2\si+2}}^{2\si+2}
+\int_{\R^d}V(x) |u(t,x)|^2dx \)= 0.
\end{align*}
The question we ask is: What remains when we
do not assume $V|a_0|^2\in L^1(\R^d)$? Roughly speaking, the 
local existence results remain when $V$ is sub-linear, but fail when
$V$ is super-linear (we prove the latter under slightly more
restrictive assumptions on $V$, see Th.~\ref{theo:superlinear}). Note
that in the above example, if we assume 
$0<\si<2/d$, then one can prove the existence of a global solution,
with an $L^2$ regularity, as in \cite{TsutsumiL2}. Our goal is to
understand better the relevance of Sobolev spaces with positive index,
when no extra decay of the initial datum is assumed. 
\smallbreak

We recall a particular case of \cite[Lemma~1]{CaBKW}:
\begin{lemma}\label{lem:hj}
  There exist $T>0$ and a unique
  solution $\phi_{\rm eik}\in C^\infty([0,T]\times\R^n)$ to:
  \begin{equation}
    \label{eq:eik}
    \partial_t \phi_{\rm eik} +\frac{1}{2}|\nabla_x \phi_{\rm eik}|^2
    +V=0\quad ;\quad 
    \phi_{{\rm eik} \mid t=0}=0\, .
  \end{equation}
This solution is sub-quadratic: $\partial_x^\alpha \phi_{\rm eik} \in
  L^\infty([0,T]\times\R^n)$ as soon as $|\alpha|\ge 2$. 
\end{lemma}
\begin{example}
  If $V(x) =\frac{1}{2}\sum_{j=1}^d \omega_j^2 x_j^2$ with $\om_j\ge 0$, then 
  \begin{equation*}
   \phi_{\rm eik}(t,x) = -\sum_{j=1}^d \frac{\omega_j}{2}x_j^2 \tan
   \(\om_j t\). 
  \end{equation*}
This shows that in general, the above result is really local in time,
due to the formation of caustics. 
\end{example}
\begin{example}
  If $V(x)= \< x\>^a$, with $0<a\le 2$, then we can see that caustics
  appear in finite time even if the potential $V$ is sub-linear.  
\end{example}
\begin{proposition}\label{prop:local}
  Let $d\ge 1$. \\
$(1)$ If $f\equiv 0$ (linear equation), assume that $a_0\in H^s(\R^d)$
for some $s\ge 0$. Then \eqref{eq:nlsp} has a unique solution $u$ such
that $u\cdot e^{-i\phi_{\rm eik}}\in C([0,T];H^s)$, where $\phi_{\rm eik}$
and $T$ are given by Lemma~\ref{lem:hj}.\\
$(2)$ For the nonlinear equation, assume that $f$ is smooth, $f\in
C^{\infty}(\R_+;\C)$, and that $a_0\in H^s(\R^d)$ for some $s>d/2$. 
Then \eqref{eq:nlsp} has a unique solution $u$ such
that $u\cdot e^{-i\phi_{\rm eik}}\in C([0,T];H^s)$, where $\phi_{\rm eik}$
and $T$ are given by Lemma~\ref{lem:hj}.
\end{proposition}
\begin{proposition}\label{prop:sublinear}
  Let $d\ge 1$,  $a_0\in H^1(\R^d)$, and assume that $V$ is sub-linear
  and that the nonlinearity
  $f$ is of the form 
  \begin{equation*}
    f(|u|^2)=\mu |u|^{2\si},\quad \text{with }\mu\in \R, \
    \si>0, \text{ and }\si<\frac{2}{d-2}\text{ if }d\ge 3.
  \end{equation*}
Then there exists $\tau=\tau(d,\|a_0\|_{H^1}, \mu,\si)>0$ such that
\eqref{eq:nlsp} has a unique solution $u\in C([0,\tau];H^1)\cap
L^{\frac{4\si+4}{d\si}}\([0,\tau]; 
W^{1,2\si+2}\)$. \\
If moreover $\si<2/d$ or  $\mu\ge 0$, then this
solution is global in time:
\begin{equation*}
 u\in C(\R_+;H^1)\cap
L^{\frac{4\si+4}{d\si}}_{\rm loc}\(\R_+; 
W^{1,2\si+2}\). 
\end{equation*}
\end{proposition}
\begin{remark*}
  Even the local result is not a consequence of
  Proposition~\ref{prop:local}: the regularity required on
  the initial data is not the same. The reason is that
  Proposition~\ref{prop:local} is established without dispersive or
  Strichartz estimates, while the local existence result in
  Proposition~\ref{prop:sublinear} is proven thanks to (local 
  in time) Strichartz estimates.
\end{remark*}

We also discuss the local Cauchy problem in $H^s(\R^d)$, $s>0$, in
Section~\ref{sec:Hs}. The main point
consists in showing that in the presence of a sub-linear potential,
local Strichartz estimates are available in Sobolev and Besov
spaces. We prove:
\begin{proposition}\label{prop:Hs}
Let $V$ be sub-linear, $0< s<d/2$ and $0<\si\le \frac{2}{d-2s}$. If
$\si$ is not an integer, assume that  $[s]<2\si$. Then 
there exist $T>0$ and a unique solution 
$$u\in C([0,T];H^s)\cap L^\gamma ([0,T];B_{\rho,2}^s)$$
to \eqref{eq:duhamel}, where
\begin{equation*}
\rho = \frac{2\si +2}{1+\frac{2\si s}{d}}\quad ;\quad \gamma =
\frac{4\si +4}{\si(d-2s)}\cdot 
\end{equation*}
\end{proposition}
We now come to the non-existence result:
\begin{theorem}\label{theo:superlinear}
  Let $d\ge 1$, and 
  $f$ be smooth, $f\in C^\infty(\R_+;\R)$. Assume that $V$ is
  super-linear, and that there exist $0<k(\le 1)$ and $C>0$ such that
  \begin{equation*}
 |\nabla V(x)|\le C \<x\>^k,\quad \forall x\in \R^d,   
  \end{equation*}
and $\omega,\omega'\in \Sph^{d-1}$ such that 
  \begin{equation}\label{eq:growth}
    |\omega\cdot\nabla V(x)|\ge c|\omega'\cdot x|^k\text{ as }|x|\to
     \infty,\text{ for some }c>0. 
  \end{equation}
Then there exists $a_0\in
  H^\infty(\R^d)$ such that
for arbitrarily small $t>0$ and all $s>0$, the solution $u(t,\cdot)$ to
  \eqref{eq:nlsp} provided by Proposition~\ref{prop:local} fails to be
  in $H^{s}(\R^d)$.  
\end{theorem}
\begin{example}
  As a potential $V$, we may consider any non-trivial quadratic form,  
or  $V(x)=\pm \<x'\>^a$, with $1<a\le 2$, for some decomposition
  $x=(x',x'')$.  
\end{example}
\begin{remark*}
  Note that no assumption is made on the growth of the nonlinearity
  at infinity: the above result reveals a \emph{geometric phenomenon},
  and not an ill-posedness result like for super-critical
  nonlinearities without a potential (\cite{BGTENS,CaARMA,CCT2,Lebeau01}). 
\end{remark*}
In Section~\ref{sec:prelim}, we outline the proof of
Proposition~\ref{prop:local}, which is a particular case of
\cite[Proposition~3]{CaBKW}. We establish Proposition~\ref{prop:sublinear}
in Section~\ref{sec:sublinear}. We extend the local theory to all the
spaces $H^s(\R^d)$ for $s>0$ in Section~\ref{sec:Hs}, where we prove
Proposition~\ref{prop:Hs}.  Finally, 
Theorem~\ref{theo:superlinear} is proved in Section~\ref{sec:general}.

\section{Preliminary remarks}
\label{sec:prelim}

In this section, we outline the proof of Proposition~\ref{prop:local},
which is a straightforward consequence of the analysis in
\cite{CaBKW}, with the choice $\varepsilon =1$. This will also guide
us for the proof of Theorem~\ref{theo:superlinear}. 
\smallbreak

First, Lemma~\ref{lem:hj} is a straightforward consequence of the
local Hamilton-Jacobi theory, Gronwall lemma, and a global inversion
theorem, which can be found for instance in
\cite[Th.~1.22]{SchwartzBook} or \cite[Prop.~A.7.1]{DG}. To prepare
the proof of Theorem~\ref{theo:superlinear}, we recall some
details.  Let $x(t,y)$ and $\xi(t,y)$ solve
\begin{equation}
  \label{eq:hamilton}
\left\{
  \begin{aligned}
   &\partial_t x(t,y) = \xi \left(t,y\right)\quad ;\quad x(0,y)=y,\\ 
   &\partial_t \xi(t,y) = -\nabla_x V\left(x(t,y)\right)\quad ;\quad
   \xi(0,y)=0.
  \end{aligned}
\right.
\end{equation}
The local Hamilton-Jacobi theory provides a solution to \eqref{eq:eik}
in the neighborhood of every point  where $y\mapsto x(t,y)$ is
invertible. The theory is global in space (not in time, in general)
thanks to the global inversion 
theorem mentioned above, and to Gronwall lemma. The gradient of
$\phi_{\rm eik}$ is given by 
\begin{equation}\label{eq:gradeik}
  \nabla_x \phi_{\rm eik}(t,x) =\xi(t,y(t,x)), 
\end{equation}
where $y(t,x)$ is the inverse mapping of $y\mapsto x(t,y)$. 
Introduce the Jacobi determinant
\begin{equation}\label{eq:defjacobi}
  J_t(y) ={\rm det}\nabla_y x(t,y). 
\end{equation}
The global inversion theorem can be applied since there
exists $C>0$ such that  
\begin{equation}\label{eq:detjacobi}
  C^{-1}\le J_t(y) \le C,\quad \forall (t,y)\in [0,T]\times \R^d.
\end{equation}
The change of unknown function $u(t,x)=
a(t,x)e^{i\phi_{\rm eik}(t,x)}$ turns \eqref{eq:nlsp} into the
equivalent Cauchy problem:
\begin{equation}
  \label{eq:hyperbolic}
    \d_t a + \nabla \phi_{\rm eik}\cdot \nabla a+\frac{1}{2}a\Delta
    \phi_{\rm eik} = \frac{i}{2}\Delta a -i f\(|a|^2\)a\quad ;\quad
a_{\mid t=0}=a_0. 
\end{equation}
The major difference with \eqref{eq:nlsp} is that the potential $V$ is
no longer present in the equation. 
The idea is to view the left hand side as a transport operator
with velocity $\nabla \phi_{\rm eik}$ and a renormalization factor
along the characteristics, $\frac{1}{2}a\Delta  \phi_{\rm eik}$. We
can then reduce the problem of 
existence of solutions of \eqref{eq:hyperbolic}, to the existence of
\emph{a priori} estimates, thanks to a mollification procedure. Since
we seek $a\in C([0,T];H^s)$, we note that the term $i\Delta$ on the
right hand side is skew-symmetric, and has no contribution in the energy
estimates. To take advantage of this property, we do not rewrite
\eqref{eq:hyperbolic} along the characteristics, but notice that from
Lemma~\ref{lem:hj}, $\|a\Delta
    \phi_{\rm eik}\|_{L^\infty_tH^s}\lesssim \|a\|_{L^\infty_tH^s}
    $. For the convective 
    term, we use Lemma~\ref{lem:hj}, and an integration by parts: if
$\alpha \in \N^d$ is such that $|\alpha|\le s$, we write
\begin{align*}
  \RE \int \d^\alpha_x \bar a \d^\alpha_x \(\nabla \phi_{\rm eik}\cdot
  \nabla a\) dx =& \RE \int \d^\alpha_x \bar a \(\nabla
  \phi_{\rm eik}\cdot 
  \nabla \d^\alpha_x a\) dx \\
&+\sum_{|\beta|\ge 1}c_{\alpha,\beta} \RE
  \int \d^\alpha_x \bar a  \nabla 
  \d^\beta_x\phi_{\rm eik}\cdot 
   \nabla \d^{\alpha-\beta}_x a dx \\
=&\frac{1}{2}\int  \nabla
  \phi_{\rm eik}\cdot 
  \nabla |\d^\alpha_x a|^2 dx +\O \(\|a\|^2_{L^\infty_t H^s}\)\\
=&\O
  \(\|a\|^2_{L^\infty_t H^s}\) .
\end{align*}
If $s$ is not an integer, we can use
interpolation. Proposition~\ref{prop:local} follows easily, since
$s>d/2$ ensures that $H^s(\R^d)$ is an algebra. 
\begin{remark}\label{rk:op}
  Let $I\subset [0,T]$ be a compact time interval. The approach of \cite{CaBKW}
  recalled above shows that the map  
$F\mapsto {\tt a}$, where
\begin{equation*}
  \d_t {\tt a} + \nabla \phi_{\rm eik}\cdot \nabla {\tt
    a}+\frac{1}{2}{\tt a}\Delta 
    \phi_{\rm eik} = \frac{i}{2}\Delta {\tt a} +F \quad ;\quad
{\tt a}_{\mid t=0}=0,
\end{equation*}
sends $L^1(I;L^2)$ to $C\cap L^\infty(I;L^2)$ continuously:
\begin{equation*}
  \|{\tt a}\|_{L^\infty(I;L^2)}\le C \|F\|_{L^1(I;L^2)},
\end{equation*}
where $C$ depends only on $d$ and $\|\nabla^2 \phi_{\rm
  eik}\|_{L^\infty(I;L^2)}$. 
\end{remark}

\section{Sub-linear potentials}
\label{sec:sublinear}

\subsection{Local $H^1$ theory}
\label{sec:local}
To prove the first part of Proposition~\ref{prop:sublinear}, the idea
is to keep the same proof as without potential. 
The gradient does not commute with $H$, but we have:
\begin{equation*}
  \(i\d_t +\frac{1}{2}\Delta\)\nabla u = V(x)\nabla u + u\nabla V(x) +
  \mu \nabla \(|u|^{2\si}u\).  
\end{equation*}
The new term is $u\nabla V(x)$, that is, $u$ multiplied by a bounded
term. Recall that $U(t)= e^{-itH}$. We show that for $\tau>0$
sufficiently small, there exists $u$ such that:
\begin{equation}
  \label{eq:duhamel}
  u(t) = U(t)a_0 -i\mu \int_0^t U(t-s)\(|u|^{2\si}u\)(s)ds=:\Phi(u)(t).
\end{equation}
We see that
\begin{equation}\label{eq:commute1}
\begin{aligned}
  \nabla \Phi(u)(t) &= U(t)\nabla a_0 -i\mu
  \int_0^t U(t-s)\nabla\(|u|^{2\si}u\)(s)ds\\
&-i  \int_0^tU(t-s)\(\Phi(u)(s)\nabla V\)ds
. 
\end{aligned}
\end{equation}
Recall that $(q,r)$ is Schr\"odinger-admissible in $\R^d$ if 
\begin{equation*}
  \frac{2}{q}+\frac{d}{r}=\frac{d}{2}\virgp\quad 2\le r\le
  \frac{2d}{d-2},\ (q,r)\not =(2,\infty).
\end{equation*}
It follows from \cite{Fujiwara} that Strichartz estimates are
available for $U(t)$ (see e.g. \cite{CazCourant}): for all
admissible pairs $(q,r)$, $(q_1,r_1)$ and $(q_2,r_2)$, there exist
$C_r$ and $C_{r_1,r_2}$ such that for any compact interval $I$ and any
$\varphi\in L^2(\R^d)$, $F\in L^{q_2}(I;L^{r_2}(\R^d))$, 
\begin{equation}
  \label{eq:strichartz}
  \begin{aligned}
    \left\| U(\cdot)\varphi\right\|_{L^q(I;L^r)}&\le C_r
    (1+|I|)^{\frac{1}{q}} \|\varphi\|_{L^2},\\
\left\| \int_0^t U(t-s)F(s)ds\right\|_{L^{q_1}(I;L^{p_1})}&\le
    C_{r_1,r_2} (1+|I|)^{\frac{1}{q_1}} 
\|F\|_{L^{q_2'}(I;L^{p_2'})}, 
  \end{aligned}
\end{equation}
where $r'$ stands for the H\"older conjugate exponent of $r$. 
Note that the powers of $|I|$ on the right hand sides are sharp in
general, for $H$ may have eigenvalues. 
For $(q,r)$ an admissible pair, define 
\begin{align*}
Y_{r,{\rm loc}}(I) := \left\{ u\in C(I;H^1); \ A u \in 
L^q_{\rm loc}(I;L^r)\cap L^\infty_{loc}(I;L^2)\ \forall A\in 
\{ {\rm Id},\nabla\}\right\}.
\end{align*}
Introduce the following Lebesgue
exponents:  
\begin{equation}\label{eq:exponents}
r=2\si +2\quad ;\quad
q=\frac{4\si+4}{d\si}\quad ;\quad k=
\frac{2\si(2\si+2)}{2-(d-2)\si}\,\cdot
\end{equation}
Then $(q,r)$ is the (admissible) pair of the proposition, and 
\begin{equation*}
\frac{1}{r'}=\frac{2\si}{r}+\frac{1}{r}\quad ; \quad
\frac{1}{q'}=\frac{2\si}{k}+\frac{1}{q}\, \cdot
\end{equation*}
For $\tau >0$ and any pair $(a,b)$, we
use the notation
\[\| f\|_{L^a_\tau L^b}= \|f\|_{L^a([0,\tau];L^b)} .\]
We first prove that there exists $\tau>0$ such
that the set 
\begin{align*}
 X_\tau := \big\{ & u \in Y_{r,{\rm loc}}([0,\tau]); \
 \|u\|_{L^\infty_\tau L^2}\le 2\|a_0\|_{L^2},\ 
\|u\|_{L^q_\tau L^r}\le 2C_{r}\|a_0\|_{L^2} \ ,\\
&\|\nabla
 u\|_{L^\infty_\tau L^2}\le  
2\|\nabla a_0\|_{L^2}\ ,\quad \|\nabla u\|_{L^q_\tau L^r}\le
 2 C_{r}\|\nabla a_0\|_{L^2}\big\}   
\end{align*}
is stable under the map $\Phi$, where $C_r$ is the constant of the
homogeneous Strichartz inequality \eqref{eq:strichartz}.
Then choosing $\tau$
even smaller, $\Phi$ is a contraction on $L^q([0,\tau];L^r)$.

Let $u \in X_\tau$. For $\tau \le 1$,
\eqref{eq:strichartz} yields:
\begin{align*}
  \|\Phi(u)\|_{L_\tau^\infty L^2}&\le \| U(t)a_0\|_{L^2} +
  C_{2,r }|\mu | \left\||u|^{2\sigma}u
  \right\|_{L^{q'}_\tau L^{r'}}\\
&\le \|a_0\|_{L^2} + C  \left\|u
  \right\|^{2\sigma}_{L^k_\tau L^r }\left\| u
  \right\|_{L^{q}_\tau L^{r} }\ .
\end{align*}
Sobolev embedding yields:
\begin{equation*}
 \left\|u
  \right\|_{L^k_\tau L^r} \le C \tau^{\frac{1}{k}}\|u\|_{L^\infty_\tau
  H^1}.
\end{equation*}
It follows that 
\begin{equation*}
  \|\Phi(u)\|_{L^\infty_\tau L^2} \le \|a_0\|_{L^2}
 + C\tau^{\frac{2\si}{k}}\|a_0\|_{H^1}^{2\si+1}.
\end{equation*}
The same computations yield:
\begin{equation*}
  \|\Phi(u)\|_{L^q_\tau L^r} \le C_r \|a_0\|_{L^2}
 + \widetilde C\tau^{\frac{2\si}{k}}\|a_0\|_{H^1}^{2\si+1}.
\end{equation*}
Similarly,
\begin{align*}
  \|\nabla \Phi(u)\|_{L_\tau^\infty L^2}\le & \, \| \nabla a_0\|_{L^2} +
  C_{2,r }|\mu | \left\|\nabla\(|u|^{2\sigma}u\)
  \right\|_{L^{q'}_\tau L^{r'}}
+ C_{2,2}\|\nabla
  V\|_{L^\infty}\|\Phi(u)\|_{L^1_\tau L^2}\\
\le &\, \|\nabla a_0\|_{L^2}  + C
  \left\|u 
  \right\|^{2\sigma}_{L^k_\tau L^r }\left\|\nabla u
  \right\|_{L^{q}_\tau L^{r} }+\\
&+\tau C_{2,2} \|\nabla V\|_{L^\infty} \(\|a_0\|_{L^2}
 + C\tau^{\frac{2\si}{k}}\|a_0\|_{H^1}^{2\si+1}\),
\end{align*}
and
\begin{align*}
  \|\nabla \Phi(u)\|_{L_\tau^q L^r}\le &\, C_r \|\nabla a_0\|_{L^2}  +
  \widetilde C
  \left\|u 
  \right\|^{2\sigma}_{L^k_\tau L^r }\left\|\nabla u
  \right\|_{L^{q}_\tau L^{r} }+\\
&+ \tau  C_{r,2} \|\nabla V\|_{L^\infty} \( \|a_0\|_{L^2}
 + C\tau^{\frac{2\si}{k}}\|a_0\|_{H^1}^{2\si+1}\).
\end{align*}
Therefore $\Phi$ leaves $X_\tau$ stable for 
\begin{equation}\label{eq:condtau}
\tau \|\nabla V\|_{L^\infty}
  \|a_0\|_{L^2}+\tau^{\frac{1}{k}}\|u\|_{L^\infty_\tau  
  H^1}\ll 1.  
\end{equation}
To complete the proof of the first part of the proposition, it is
enough to prove contraction for small $\tau$ in 
the weaker metric $L^q([0,\tau];L^r)$. 
We have: 
\begin{equation*}
\begin{aligned}
    \big\| \Phi(u_2)-\Phi(u_1)
    \big\|_{L^q_\tau L^r }  & \le
C\left\| \left(|u_2|^{2\sigma} u_2 - |u_1|^{2\sigma}
    u_1\right) \right\|_{L^{q'}_\tau L^{r'} }\\
&\le C\left(
    \|u_1\|^{2\sigma}_{L^k_\tau L^r }+
\|u_2\|^{2\sigma}_{L^k_\tau L^r }\right) 
\|u_2-u_1\|_{L^q_\tau L^r }\, .
\end{aligned}
\end{equation*}
As above, we have the estimate
\[
\|u_j\|^{2\sigma}_{L^k_\tau L^r }\le C
\tau^{\frac{1}{k}}\|u_j\|_{L^\infty_\tau 
  H^1}. 
\]
Therefore,  contraction follows for $\tau$ sufficiently small,
according to \eqref{eq:condtau}.  

\subsection{Global existence in $H^1$}
\label{sec:globalH1}

If $V$ is sub-linear and unbounded, then the energy
\begin{equation*}
 E= \frac{1}{2}\|\nabla
u(t)\|_{L^2}^2+\frac{\mu}{\si+1} \|u(t)\|_{L^{2\si+2}}^{2\si+2}
+\int_{\R^d}V(x) |u(t,x)|^2dx
\end{equation*}
may not be defined initially, if we simply require $a_0\in
H^1(\R^d)$. To complete the proof of 
Proposition~\ref{prop:sublinear}, the idea is to notice that the time
derivative of the ``bad'' term in the energy is controlled by the
$H^1$ norm of the solution. We present the computations at a formal
level only, and refer to \cite{CazCourant} for a justification
method which uses the multiplication by Gaussians. We have
\begin{align*}
  \frac{d}{dt}\int_{\R^d}V(x) |u(t,x)|^2dx &= 2\RE \int_{\R^d}V(x)
  \overline{u} \d_t udx =2\IM\int_{\R^d}V(x)
  \overline{u} \(i\d_t u\)dx\\
&= -\IM\int_{\R^d}V(x)
  \overline{u} \Delta udx  = \IM\int_{\R^d}\overline{u}\nabla V(x)
  \cdot \nabla  udx.
\end{align*}
We infer, thanks to the conservation of mass:
\begin{align*}
  \frac{1}{2}\|\nabla
u(t)\|_{L^2}^2+\frac{\mu}{\si+1} \|u(t)\|_{L^{2\si+2}}^{2\si+2} \le &
\frac{1}{2}\|\nabla 
a_0\|_{L^2}^2+\frac{\mu}{\si+1} \|a_0\|_{L^{2\si+2}}^{2\si+2}\\
&+ 
\|\nabla V\|_{L^\infty}\|a_0\|_{L^2}\int_0^t \|\nabla u(s)\|_{L^2}ds.   
\end{align*}
When $\mu \ge 0$, this yields the estimate
\begin{equation*}
  \|\nabla
u(t)\|_{L^2}^2 \lesssim 1 + \int_0^t \|\nabla u(s)\|_{L^2}ds, 
\end{equation*}
hence $\|\nabla u(t)\|_{L^2}$ grows at most exponentially. 
\smallbreak

If $\si<2/d$ and $\mu <0$, Gagliardo--Nirenberg inequality and the
conservation of mass yield:
\begin{align*}
  \|\nabla
u(t)\|_{L^2}^2 \lesssim 1 + \|\nabla
u(t)\|_{L^2}^{d\si} + \int_0^t \|\nabla u(s)\|_{L^2}ds.
\end{align*}
Using Young inequality
\begin{equation*}
  \|\nabla
u(t)\|_{L^2}^{d\si} \le C_\epsilon + \epsilon \|\nabla
u(t)\|_{L^2}^{2}, 
\end{equation*}
and choosing $\epsilon >0$ sufficiently small, we conclude as
before. This completes the proof of Proposition~\ref{prop:sublinear}. 

\section{On the local Cauchy problem in $H^s$: proof of
  Proposition~\ref{prop:Hs}} 
\label{sec:Hs}
When $a_0\in H^s(\R^d)$ with $s>0$ not necessarily equal to one, and
$V$ is sub-linear, it is
still possible to establish a local in time theory. 
Without potential, 
$V\equiv 0$, Proposition~\ref{prop:Hs} was proved by T.~Cazenave and
F.~Weissler  
\cite[Theorem~1.1, (i)--(ii)]{CW90}.  As in this paper, we shall not
define Besov spaces by 
using a dyadic decomposition, but rather use their characterization in
terms of interpolation between Sobolev spaces. 
We first recall the argument when $V\equiv 0$, and then show how
it can be adapted to infer Proposition~\ref{prop:Hs}.

\subsection{Proof when $V\equiv 0$}

The idea is to apply a fixed point argument, as in
Section~\ref{sec:local}. However, when $s<d/2$ is not an integer, it
becomes delicate to estimate the $H^s$ norm of the
nonlinearity. This is why in \cite{CW90}, the authors work in Besov
spaces. 
When $s$ is an integer, the above result can be refined. We shall not
recall this aspect more precisely, and simply refer to
\cite{CW90}. The proof proceeds in
three steps. The authors first establish Strichartz estimates for the
free group $e^{i\frac{t}{2}\Delta}$ in (homogeneous) Besov
spaces \cite[Th.~2.2]{CW90}. Next, they prove 
estimates for the nonlinear term, in homogeneous Besov spaces as
well \cite[Th.~3.1]{CW90}. Finally, these tools, along with Strichartz
estimates, make it 
possible to apply a fixed point argument, and prove
Proposition~\ref{prop:Hs} when $V\equiv 0$.
\smallbreak

Denote
\begin{equation*}
{\mathcal I}(t)F := \int_0^t U(t-s)F(s)ds.
\end{equation*}
The first step yields, for $s>0$, and $(q,r)$, $(q_j,r_j)$ admissible
pairs: 
\begin{equation}
  \label{eq:strichartzCW90}
  \begin{aligned}
    \left\| U(\cdot)\varphi\right\|_{L^q(\R_+;\dot B_{r,2}^s)}&\le C_r
    \|\varphi\|_{\dot H^s},\\
\left\| {\mathcal I}(\cdot)F\right\|_{L^{q_1}(I;\dot B_{r_1,2}^s)}&\le
    C_{r_1,r_2} 
\|F\|_{L^{q_2'}(I;\dot B_{r_2,2}^s)}, 
  \end{aligned}
\end{equation}
where $C_{r_1,r_2}$ does not depend on the time interval $I$. 
Next, under the assumptions of Proposition~\ref{prop:Hs}, we have
\begin{equation}\label{eq:estnlCW90}
\left\| |u|^{2\si}u\right\|_{\dot B_{\rho',2}^s}\lesssim \left\|
u\right\|_{\dot B_{\rho,2}^s}^{2\si+1}.
\end{equation}
Proposition~\ref{prop:Hs} follows from \eqref{eq:strichartzCW90},
\eqref{eq:estnlCW90}, H\"older's inequality and a fixed point
argument. 
\begin{remark}
Note that \eqref{eq:strichartzCW90} and 
\eqref{eq:estnlCW90} still hold if we replace homogeneous Besov spaces
with inhomogeneous ones. This remark simplifies the generalization to
the case when $V$ is sub-linear.  
\end{remark}
\subsection{Strichartz estimates in Besov spaces with a sub-linear
  potential} 
We show that when $V$ is sub-linear, \eqref{eq:strichartzCW90} still
holds, up to two modifications:
\begin{itemize}
\item The Strichartz inequalities hold on \emph{finite} time
  intervals only.
\item We replace the homogeneous Besov spaces
with inhomogeneous ones.
\end{itemize}
The first point is unavoidable, as recalled in
Section~\ref{sec:local}. Since we shall prove a local in time result,
in the rest of this section we consider time intervals of
length at most one. The second point is here to consider
pseudo-differential operators with smooth symbols which do not contain
$x$-variable. 
\smallbreak

If $P=P(D)$ is a pseudo-differential operator with smooth symbol, we
have:
\begin{align*}
[P,U(t)]\varphi &= -i\int_0^t U(t-s)[P,V]U(s)\varphi ds= -i {\mathcal
  I}(t) \( [P,V]U(\cdot)\varphi\),\\
[P,{\mathcal I}(t)]F &= -i\int_0^t U(t-s)[P,V]{\mathcal I}(s)F ds=-i {\mathcal
  I}(t) \( [P,V]{\mathcal I}(\cdot)\varphi\).
\end{align*}
First, assume $0<s<1$. For $I$ a time interval with $|I|\le 1$,
\eqref{eq:strichartz} yields:
\begin{align*}
\left\|P U(t)\varphi \right\|_{L^q(I;L^r)}&\le \left\| U(t)P\varphi
\right\|_{L^q(I;L^r)} + \left\|{\mathcal
  I}(t)\( [P,V]U(\cdot)\varphi\) \right\|_{L^q(I;L^r)}\\
&\lesssim \|P\varphi\|_{L^2} + \left\| [P,V]U(\cdot)\varphi
\right\|_{L^1(I;L^2)}\\
&\lesssim  \|P\varphi\|_{L^2} + \left\|
[P,V]U(\cdot)\varphi \right\|_{L^\infty(I;L^2)}.
\end{align*}
Similarly,
\begin{equation*}
\left\|P {\mathcal I}(t)F\right\|_{L^{q_1}(I;L^{r_1})}\lesssim
\|P F\|_{L^{q'_2}(I;L^{r'_2})}+ \left\|
[P,V]{\mathcal I}(\cdot)F \right\|_{L^\infty(I;L^2)}. 
\end{equation*}
For $s>0$, let $P_s =(I-\Delta)^{s/2}$. By
\cite[Th.~2]{CoifmanMeyerAIF} (see also \cite[\S~3.6]{Taylor91}), we
know that if in addition $s\le 1$, then $[P_s,V]$ is bounded from
$L^2$ to $L^2$, with norm  controlled by $C\|\nabla V\|_{L^\infty}$
for some universal constant $C$. We infer, when $s\le 1$, 
\begin{equation*}
\begin{aligned}
\left\|P_s U(t)\varphi \right\|_{L^q(I;L^r)}&\lesssim  \|P_s\varphi\|_{L^2}
+ \left\| U(\cdot)\varphi \right\|_{L^\infty(I;L^2)}\\
&\lesssim \|P_s\varphi\|_{L^2} 
+ \|\varphi\|_{L^2}\lesssim \|P_s\varphi\|_{L^2} ,
\end{aligned}
\end{equation*}
where we have used Strichartz estimates
\eqref{eq:strichartz}. This means:
\begin{equation}\label{eq:striHshomo}
\left\|U(\cdot)\varphi \right\|_{L^q(I;W^{s,r})}\lesssim \|\varphi\|_{H^s} 
\end{equation}
Similarly, when $s\le 1$, 
\begin{equation}\label{eq:striHsinhomo}
\left\|{\mathcal I}(\cdot)F\right\|_{L^{q_1}(I;W^{s,r_1})}\lesssim
\| F\|_{L^{q'_2}(I;W^{s,r'_2})}. 
\end{equation}
For $s>1$, replace $P_s$ with the family $(P_{s-m}\circ\d^\alpha
)_{|\alpha|\le m}$, where $m=[s]$. Reasoning as above, we see that 
since $\d^\alpha V\in L^\infty(\R^d)$ for all $|\alpha|\ge 1$,
\eqref{eq:striHshomo} and \eqref{eq:striHsinhomo} hold for all
$s>0$. 
\smallbreak

Interpolating (as in \cite{CW90}, up to replacing homogeneous spaces
by their inhomogeneous counterparts), we conclude: 
\begin{equation}
  \label{eq:strichartzBesov}
  \begin{aligned}
    \left\| U(\cdot)\varphi\right\|_{L^q(I; B_{r,2}^s)}&\le C_r
    \|\varphi\|_{H^s},\\
\left\| {\mathcal I}(\cdot)F\right\|_{L^{q_1}(I; B_{r_1,2}^s)}&\le
    C_{r_1,r_2} 
\|F\|_{L^{q_2'}(I; B_{r_2,2}^s)}, 
  \end{aligned}
\end{equation}
where the constants $C_r$ and $C_{r_1,r_2}$ do not depend on $I$,
provided that $|I|\le 1$. 
\bigbreak

\noindent{\bf Conclusion.} Since \eqref{eq:estnlCW90} holds with
homogeneous Besov spaces replaced by their inhomogeneous counterparts,
the fixed point argument used in \cite{CW90} can be transported
here. This completes the proof of Proposition~\ref{prop:Hs}.

\section{Loss of Sobolev regularity: proof of
  Theorem~\ref{theo:superlinear}}
\label{sec:general}

\subsection{A decomposition suggested by geometric optics}
\label{sec:decomp}
The idea is to resume the approach of weakly nonlinear geometric
optics recalled in Section~\ref{sec:prelim}. We consider an
intermediary function defined by leaving out the term $i\Delta a$ in
\eqref{eq:hyperbolic}: without this term,
\eqref{eq:hyperbolic} is  an ordinary differential equation
along the characteristics of the transport operator with velocity
$\nabla \phi_{\rm eik}$ (\emph{i.e.} the bicharacteristics associated to $H$). 
\smallbreak

Recall that $a$ solves
\eqref{eq:hyperbolic}, and define $b$ as the solution on $[0,T]$ to:
\begin{equation}
  \label{eq:bgen}
    \d_t b + \nabla \phi_{\rm eik}\cdot \nabla b+\frac{1}{2}b\Delta
    \phi_{\rm eik} = -i f\(|b|^2\)b\quad ;\quad
b_{\mid t=0}=a_0. 
\end{equation}
To see that $b$ solves an ordinary differential equation along the
rays of geometric optics (the projections of the Hamilton flow
\eqref{eq:hamilton} on the physical space), introduce
\begin{equation*}
  \beta(t,y) =b\(t,x(t,y)\)\sqrt{J_t(y)}, 
\end{equation*}
where $x(t,y)$ is given by \eqref{eq:hamilton} and the Jacobi
determinant is defined by \eqref{eq:defjacobi}. This change of unknown
function makes sense for $t\in [0,T]$, where $y\mapsto x(t,y)$ is a
global diffeomorphism. Then \eqref{eq:bgen} is equivalent to 
\begin{equation}\label{eq:beta}
  \d_t \beta(t,y) = -i f\(J_t(y)^{-1}|\beta(t,y)|^2\)\beta(t,y)\quad
  ;\quad \beta(0,y)=a_0(y).  
\end{equation}
Since in Theorem~\ref{theo:superlinear}, we assume that $f$ is
real-valued, we note that 
\begin{equation*}
  \d_t|\beta|^2 =0,
\end{equation*}
so that \eqref{eq:beta} is just a \emph{linear} ordinary differential
equation:
\begin{equation*}
  \beta(t,y)= a_0(y)\exp\(-i \int_0^t f\(J_s(y)^{-1}|a_0(y)|^2\)ds\).
\end{equation*}
We infer
\begin{equation*}
  b(t,x) = \frac{1}{\sqrt{J_t(y(t,x))}} a_0\( y(t,x)\) \exp\(-i \int_0^t
  f\(J_s\(y(t,x)\)^{-1}\left|a_0\(y(t,x)\)\right|^2\)ds\) .
\end{equation*}
The main observation is that \eqref{eq:detjacobi} implies that $b\in
C([0,T];H^s(\R^d))$ for all $s\ge 0$. Let $r=a-b$: for every $t\in
[0,T]$, $r(t,\cdot)\in H^\infty(\R^d)$. For $1\le j\le d$, $x_j r$
solves:
\begin{align*}
  \d_t(x_j r) +\nabla \phi_{\rm eik}\cdot \nabla (x_j
    r)+\frac{1}{2}x_j r \Delta
    \phi_{\rm eik} &= \frac{i}{2}\Delta (x_j r)+r\d_j \phi_{\rm
    eik}-i\d_j r+\frac{i}{2}x_j\Delta b \\
&-i   x_j\(f\(|b+r|^2\)(b+r)- f\(|b|^2\)b\),\\
x_j r_{\mid t=0}&= 0.
\end{align*}
Notice that the fundamental theorem of calculus yields:
\begin{align*}
  x_j \( f\(|a|^2\)a-f\(|b|^2\)b\) &= x_j \(
  f\(|b+r|^2\)(b+r)-f\(|b|^2\)b\)\\
&= x_j r\int_0^1 \d_z F\( b+sr\)ds + x_j \overline{r}\int_0^1
  \d_{\overline z} F\( b+sr\)ds ,
\end{align*}
where $F(z)=f(|z|^2)z$. In particular, we know that 
\begin{align*}
  \int_0^1\d_z F\( b+sr\)ds , \int_0^1
  \d_{\overline z} F\( b+sr\)ds\in C\cap L^\infty (I\times \R^d).
\end{align*}
Reasoning as in Remark~\ref{rk:op}, we see that:
\begin{equation*}
  \|x r\|_{L^\infty([0,t];L^2)}\le C\(1 + \|x\Delta
  b\|_{L^1([0,t];L^2)} \) .
\end{equation*}
We must make sure that the last term is, or can be chosen, finite. 
We shall demand $x\Delta b \in L^\infty([0,T];L^2)$. In view of
\eqref{eq:detjacobi}, 
this requirement is met as soon as $a_0\in H^\infty(\R^d)$ is such
that $x\Delta a_0,x a_0|\nabla a_0|^2\in L^2(\R^d)$. We then have:
\begin{equation}\label{eq:17h15}
\begin{aligned}
&\text{If }a_0\in H^\infty(\R^d)\text{ is such that }x\Delta a_0,x
 a_0|\nabla a_0|^2\in L^2(\R^d), \text{ then:}\\
& a= b+r,\text{ with }b,r\in C([0,T];H^s)
 \ \forall s\ge 0, \text{ and }xr \in
C([0,T];L^2).
\end{aligned}
\end{equation}

\subsection{Small time approximation of $\nabla \phi_{\rm eik}$}
\label{sec:phismalltime}
We now prove that for small times, $\nabla \phi_{\rm eik}(t,x)$
can be approximated by $-t\nabla V(x)$. 
\begin{lemma}\label{lem:dl}
Assume that there exist $0\le k\le 1$ and $C>0$ such that
  \begin{equation*}
    |\nabla V(x)|\le C\<x\>^k,\quad \forall x\in \R^d.
  \end{equation*}
Then there exist
  $T_0,C_0>0$  such that 
  \begin{equation*}
    \left| \nabla \phi_{\rm eik}(t,x) +t\nabla V(x) \right| \le
    C_0 t^2\<x\>^k,\quad \forall t\in [0,T_0].
  \end{equation*}
\end{lemma}
\begin{proof}[Proof of Lemma~\ref{lem:dl}] We infer from
    \eqref{eq:eik} and Lemma~\ref{lem:hj} that 
    \begin{equation}\label{eq:11h43}
      \left| \d_t \nabla \phi_{\rm eik}(t,x)+\nabla V(x)\right| \le
      \|\nabla^2 \phi_{\rm eik}(t)\|_{L^\infty} |\nabla \phi_{\rm
      eik}(t,x)| \lesssim |\nabla \phi_{\rm
      eik}(t,x)|. 
    \end{equation}
From \eqref{eq:hamilton} and \eqref{eq:gradeik}, we also have
\begin{align*}
  |\nabla \phi_{\rm eik}(t,x)|&=\left|\xi\(t,y(t,x)\)\right| =
   \left|\int_0^t\nabla V\(x(s,y(t,x))\)ds\right| \\
&\lesssim \int_0^t\left|\nabla V\(y(t,x)\)\right| ds
   +\int_0^t\left|x(s,y(t,x))- y(t,x)\right| ds.  
\end{align*}
We claim that 
\begin{equation}\label{eq:11h33}
  \left| x(t,y) - y\right| \lesssim t^2\<y\>^k.
\end{equation}
Indeed, we have from \eqref{eq:hamilton},
\begin{align*}
  \left|x(t,y) - y\right| &=\left|\int_0^t \d_t x(s,y)ds\right| =
  \left|\int_0^t\int_0^s \nabla V\( 
  x(s',y)\)ds'ds \right|\\
&=\left|\int_0^t (t-s')\nabla V\(
  x(s',y)\)ds'\right|\\
&= \left|\int_0^t (t-s) \nabla V\(
  y\)ds + \int_0^t (t-s)\(\nabla V\(  x(s,y)\)- \nabla V\( y\)\)ds
  \right|\\
&\lesssim t^2 \<y\>^k +  \int_0^t (t-s) \left|x(s,y) - y\right|ds,
\end{align*}
and \eqref{eq:11h33} follows from Gronwall lemma. We infer that for
$t>0$ sufficiently small,
\begin{equation*}
  \left| y(t,x) - x\right| \lesssim t^2\<x\>^k,
\end{equation*}
and therefore,
\begin{align*}
 |\nabla \phi_{\rm eik}(t,x)|&\lesssim \int_0^t\left|\nabla
   V\(y(t,x)\)\right| ds 
   +\int_0^t\left|x(s,y(t,x))- y(t,x)\right| ds\\
&\lesssim \int_0^t\left|\nabla
   V\(x\)\right| ds +\int_0^t\left|x- y(t,x)\right| ds
   +\int_0^t\left|x(s,y(t,x))- y(t,x)\right| ds\\
&\lesssim t\<x\>^k +  t^3\<x\>^k+ \int_0^t s^2\< y(t,x)\>^kds\\
&  \lesssim t\<x\>^k +  t^3\<x\>^k+ t^3 \( \<x\>^k + t^{2k}\<x\>^{2k}\).  
\end{align*}
Then \eqref{eq:11h43} yields
\begin{equation*}
  \left| \d_t \nabla \phi_{\rm eik}(t,x)+\nabla V(x)\right| \lesssim
  t\<x\>^k ,
\end{equation*}
Lemma~\ref{lem:dl} follows by integration in time.  
\end{proof}
We infer that for $t>0$ small enough,
\begin{equation}\label{eq:17h14}
  |\omega\cdot \nabla \phi_{\rm eik}(t,x)|\gtrsim t |\omega\cdot
   \nabla V(x)|.
\end{equation}
\subsection{Conclusion}
\label{sec:concl}
Consider
\begin{equation}\label{eq:defa0}
    a_0(x) =\frac{1}{\<x\>^{d/2}\log \( 2+|x|^2\)}\cdot
  \end{equation}
As is easily checked, $a_0$ meets the requirements of the first line
of \eqref{eq:17h15}. Denote
\begin{equation*}
v = b
e^{i\phi_{\rm eik}}\quad ;\quad w = re^{i\phi_{\rm eik}}. 
\end{equation*}
Obviously, $u=v+w$. 
We see from \eqref{eq:17h15} and \eqref{eq:17h14} that $v(t,\cdot)
\in L^2(\R^d)\setminus H^1(\R^d)$ 
for $t>0$ sufficiently small, under the assumptions of
Theorem~\ref{theo:superlinear}. On the other hand, $w(t,\cdot)
\in H^1(\R^d)$ for all $t\in [0,T]$, hence $u(t,\cdot)\in
L^2(\R^d)\setminus H^1(\R^d)$ 
for $0<t\ll 1$. 
\smallbreak

To complete the proof of
Theorem~\ref{theo:superlinear}, we now just have to see that the same
holds if we replace $H^1(\R^d)$ with $H^s(\R^d)$ for $0<s<1$. We use
the following characterization of $H^s(\R^d)$ (see
e.g. \cite{JYClivre}): for $\varphi\in L^2(\R^d)$ and $0<s<1$, 
\begin{equation*}
\varphi\in H^s(\R^d)\iff  \ii_{\R^d\times \R^d}
  \frac{\left|\varphi(x+y) - \varphi(x)\right|^2}{|y|^{d +2s}
  }dxdy<\infty. 
\end{equation*}
Since $w(t,\cdot)\in 
H^1$ for all $t\in [0,T]$, we shall prove that $v(t,\cdot)
\in L^2\setminus H^s$ for $t$ sufficiently small. Let $0<s<1$. We
prove that for $0<t\ll 1$, 
\begin{equation*}
I:= \int_{|y|\le 1}\int_{x\in \R^d} \frac{\left|v(t,x+y) -
  v(t,x)\right|^2}{|y|^{d +2s}}dxdy=\infty.   
\end{equation*}
To apply a fractional Leibnitz rule, write
\begin{align*}
  v(t,x+y) - v(t,x)&=\(b(t,x+y) -
  b(t,x)\)e^{i\phi_{\rm eik}(t,x+y)}\\
&\quad  + \(e^{i\phi_{\rm eik}(t,x+y)}-
  e^{i\phi_{\rm eik}(t,x)}\)b(t,x).
\end{align*}
In view of the inequality $|\alpha -\beta|^2 \ge \alpha^2/2 -\beta^2$,
we have: 
\begin{align*}
  \left|v(t,x+y) - v(t,x)\right|^2&\ge
  \frac{1}{2}\left|\(e^{i\phi_{\rm eik}(t,x+y)}- 
  e^{i\phi_{\rm eik}(t,x)}\)b(t,x) \right|^2 \\
&\quad - \left|b(t,x+y) -
  b(t,x)\right|^2.
\end{align*}
We can leave out the last term, since $b(t,\cdot)\in H^\infty$ for
$t\in [0,T]$:
\begin{equation*}
  \ii_{\R^d\times \R^d}
  \frac{\left|b(t,x+y) - b(t,x)\right|^2}{|y|^{d +2s}
  }dxdy<\infty,\quad \forall t\in[0,T]. 
\end{equation*}
We now want to prove
\begin{equation*}
  \int_{|y|\le 1}\int_{x\in \R^d} |b(t,x)|^2
\frac{\left|\sin \(\frac{\phi_{\rm eik}(t,x+y)-\phi_{\rm
        eik}(t,x)}{2}\)\right|^2}{|y|^{d +2s}}dxdy=\infty.
\end{equation*}
Lemma~\ref{lem:hj} yields: 
\begin{equation*}
 \(\d_t +\nabla \phi_{\rm eik}\cdot \nabla \)\nabla^2 \phi_{\rm
 eik}\in L^\infty \([0,T]\times\R^d\)^{d^2}\quad ;\quad \nabla^2 \phi_{{\rm
 eik} \mid t=0}=0. 
\end{equation*}
Therefore,
\begin{equation*}
\left\|  \nabla^2 \phi_{\rm
    eik}(t,\cdot)\right\|_{L^\infty(\R^d)^{d^2}}=\O(t)\quad \text{as
    }t\to 0. 
\end{equation*}
We infer:
\begin{equation*}
  \phi_{\rm eik}(t,x+y)-\phi_{\rm eik}(t,x)=y\cdot \nabla \phi_{\rm
  eik}(t,x) + \O(t|y|^2),\quad \text{uniformly for }x\in \R^d,
\end{equation*}
and 
\begin{align*}
  \sin \(\frac{\phi_{\rm eik}(t,x+y)-\phi_{\rm
        eik}(t,x)}{2}\) &= \sin \(\frac{y\cdot \nabla \phi_{\rm
  eik}(t,x)}{2}\)\cos\( \O(t|y|^2)\)\\
&\quad +\cos \(\frac{y\cdot \nabla \phi_{\rm
  eik}(t,x)}{2}\)\sin\( \O(t|y|^2)\). 
\end{align*}
The second term is $\O(t|y|^2)$. Using the estimate $|\alpha
-\beta|^2 \ge \alpha^2/2 -\beta^2$ again, we see that the integral
corresponding to the second term is finite, and can be left out. 
To prove that 
\begin{equation*}
  I'=\int_{|y|\le 1}\int_{x\in \R^d} |b(t,x)|^2
\frac{\left|\sin \(\frac{y\cdot \nabla \phi_{\rm
  eik}(t,x)}{2}\)\right|^2}{|y|^{d +2s}}dxdy=\infty \quad
\text{for }0<t\ll 1,
\end{equation*}
we can localize $y$ in a small conic
neighborhood of $\omega \R \cap \{|y|\le 1\}$: 
\begin{equation*}
{\mathcal V}_\epsilon = \{|y|\le 1\ ;\ 
|y-(y\cdot \omega) \omega|\le \epsilon |y|\}, \quad 0<\epsilon \ll 1.  
\end{equation*}
For $0<\epsilon,t\ll 1$, \eqref{eq:17h14} yields:
\begin{equation*}
 \left|\sin \(\frac{y\cdot \nabla \phi_{\rm
  eik}(t,x)}{2}\)\right|\gtrsim t|y\cdot \omega|\times
  \left|\omega \cdot \nabla 
  V(x)\right|,\quad y\in  {\mathcal V}_\epsilon.
\end{equation*}
Introduce a conic localization for $x$ close to $\omega'$, excluding
the origin:
\begin{equation*}
  {\mathcal U}_\epsilon = \{ |x|\ge 1\ ;\ |x-(x\cdot \omega')
  \omega|\le \epsilon |x|\}. 
\end{equation*}
Change the variable in the $y$-integral: for $t$ and $\epsilon$
sufficiently small, and $x\in {\mathcal U}_\epsilon$, 
set 
\begin{equation*}
  y' =  \omega\cdot \nabla \phi_{\rm
  eik}(t,x)y. 
\end{equation*}
This change of variable is admissible, from \eqref{eq:growth} and
\eqref{eq:17h14}.  
For $0<\epsilon,t\ll 1$, we have:
\begin{align*}
  I'&\ge \int_{y\in {\mathcal V}_\epsilon}\int_{x\in \R^d} |b(t,x)|^2
\frac{\left|\sin \(\frac{y\cdot \nabla \phi_{\rm
  eik}(t,x)}{2}\)\right|^2}{|y|^{d +2s}}dxdy \\
&\gtrsim \int_{x\in {\mathcal U}_\epsilon} 
|b(t,x)|^2 |\omega\cdot \nabla \phi_{\rm
  eik}(t,x)|^{2s}\(\int_{y\in |\omega\cdot \nabla \phi_{\rm
  eik}(t,x)|{\mathcal V}_\epsilon} 
\frac{dy}{|y|^{d +2s-2}}\)dx\\
&\gtrsim \int_{x\in {\mathcal U}_\epsilon} 
|b(t,x)|^2 |\omega\cdot \nabla \phi_{\rm
  eik}(t,x)|^{2s}\(\int_{y\in c t {\mathcal
    V}_\epsilon}  
\frac{dy}{|y|^{d +2s-2}}\)dx.
\end{align*}
The assumption \eqref{eq:growth}, the expression of $b$ and the
choice \eqref{eq:defa0} for $a_0$ then show that for $0<t\ll 1$,
$I=\infty$. This completes the proof of
Theorem~\ref{theo:superlinear}.  

\providecommand{\bysame}{\leavevmode\hbox to3em{\hrulefill}\thinspace}
\providecommand{\href}[2]{#2}

\end{document}